\documentclass{amsart}
\usepackage{graphics}
\usepackage{amscd}
\usepackage{amsfonts,amssymb,latexsym}  

\newenvironment{proof2}{\hspace{-\parindent}\textit{Proof}}{\hfill $\Box$}

\newcommand{\SA}{{\mathcal{A}}}

\newcommand{\SC}{{\mathcal{C}}}

\newcommand{\SF}{{\mathcal{F}}}
\newcommand{\SG}{{\mathcal{G}}}
\newcommand{\SH}{{\mathcal{H}}}

\newcommand{\SO}{{\mathcal{O}}}
\newcommand{\SP}{{\mathcal{P}}}
\newcommand{\SQ}{{\mathcal{Q}}}

\newcommand{\calS}{{\mathcal{S}}}
\newcommand{\ST}{{\mathcal{T}}}

\newcommand{\SV}{{\mathcal{V}}}

\DeclareFontFamily{OT1}{rsfs}{}
\DeclareFontShape{OT1}{rsfs}{n}{it}{<->rsfs10}{}
\DeclareMathAlphabet{\curly}{OT1}{rsfs}{n}{it}

\newcommand{\CrC}{{\curly C}}
\newcommand{\CrE}{{\curly E}}
\newcommand{\CrF}{{\curly F}}
\newcommand{\CrK}{{\curly K}}
\newcommand{\CrH}{{\curly H}}
\newcommand{\CrL}{{\curly L}}
\newcommand{\CrM}{{\curly M}}
\newcommand{\CrN}{{\curly N}}
\newcommand{\CrT}{{\curly T}}
\newcommand{\CrV}{{\curly V}}

\newcommand{\PP}{\mathbb{P}}
\newcommand{\ZZ}{\mathbb{Z}}
\newcommand{\CC}{\mathbb{C}}
\newcommand{\QQ}{\mathbb{Q}}
\newcommand{\VV}{\mathbb{V}}
\newcommand{\FM}{\mathfrak{M}}

\newcommand{\isom}{\cong}
\newcommand{\Ext}{\operatorname{Ext}}
\newcommand{\Grass}{\operatorname{Grass}}

\newcommand{\Hilb}{\operatorname{Hilb}}
\newcommand{\Hom}{\operatorname{Hom}}

\newcommand{\Pic}{\operatorname{Pic}}
\newcommand{\Sym}{\operatorname{Sym}}
\newcommand{\id}{\operatorname{id}}
\newcommand{\im}{\operatorname{im}}
\newcommand{\surj}{\twoheadrightarrow}
\newcommand{\inj}{\hookrightarrow}
\newcommand{\gr}{\operatorname{gr}}
\newcommand{\rk}{\operatorname{rk}}

\newcommand{\wt}{\widetilde}

\newcommand{\slv}{\operatorname{SL}(V)}
\newcommand{\tube}{conic bundle}
\newcommand{\tubes}{conic bundles}
\newcommand{\tq}{\tilde q}
\newcommand{\tQ}{\tilde Q}
\newcommand{\z}{(\tq,\tQ)}
\newcommand{\chim}{\chi(\SO_X(m))}

\newcommand{\ssmoduli}{\overline\FM_\tau (r,d,\CrL)}

\newcommand{\imag}{i}
\newcommand{\be}{\begin{equation}}
\newcommand{\ee}{\end{equation}}

\newcommand{\Id}{\mbox{Id}}

\newcommand{\Vol}{\mbox{Vol}}

\newcommand{\la}{\langle}
\newcommand{\ra}{\rangle}

\newtheorem{proposition}{Proposition}[section]
\newtheorem{theorem}[proposition]{Theorem}
\newtheorem{definition}[proposition]{Definition}
\newtheorem{lemma}[proposition]{Lemma}

\newtheorem{corollary}[proposition]{Corollary}
\newtheorem{remark}[proposition]{Remark}

\title[Stability of Conic Bundles]{Stability of 
Conic Bundles (with an appendix by Mundet i Riera)}
\author{T. G\'omez and I. Sols}
\date{July 2, 1999}
\thanks{Mathematical Subject Classification: Primary 14D22, 
Secondary 14D20}

\begin{document}

\begin{abstract}
Roughly speaking, a \tube\ is a surface, 
fibered over a curve, such that the fibers are
conics (not necessarily smooth). We define stability for \tubes\ and
construct a moduli space. We prove that (after fixing some invariants)
these moduli spaces are irreducible (under some conditions). Conic
bundles can be thought of as generalizations of orthogonal bundles on
curves. We show that in this particular case our definition of stability agrees
with the definition of stability for orthogonal bundles. 
Finally, in
an appendix by I. Mundet i Riera, a Hitchin-Kobayashi correspondence 
is stated for \tubes.
\end{abstract}

\maketitle

\section{Introduction}

In this paper we introduce the notion of stable \tube. This notion
appears as the stability
condition in the GIT construction of the moduli space of these
objects. 

Let $X$ be a smooth complex curve of genus $g$.
Let $r>0$ and $d$ be two integer numbers. Let $\CrL$ be a line bundle
over $X$. These data will be fixed throughout the paper.

\begin{definition}
A \tube\ on $X$ of type $(r,d,\CrL)$ is a pair $(\CrE,Q)$ where $\CrE$ is a
vector bundle on $X$ of rank $r$ and degree $d$, and $Q$ is a morphism
$$
Q: Sym^2 \CrE \to \CrL.
$$
A morphism between \tubes\ $\varphi:(\CrE,Q)\to (\CrE',Q')$ is a
morphism
$f:\CrE \to \CrE'$ such that there is a commutative diagramme
$$
\CD
\Sym^2 \CrE @>\Sym^2 f>> \Sym^2 \CrE' \\
@V{Q}VV             @V{Q'}VV \\
\CrL @>g>>   \CrL
\endCD
$$
where $g$ is a scalar multiple of identity.
\end{definition}

Then two \tubes\ $(\CrE,Q)$ and $(\CrE',Q')$ will be isomorphic when 
there is an isomorphism $\CrE\isom \CrE'$ 
that takes $Q$ into a scalar multiple of $Q'$. The name \tube\ comes
from the case $r=3$. We will be mostly interested in this case, and in
fact we will only define stability for $r \leq 3$.

If $\CrL=\SO_X$ and $Q$ gives a nondegenerate quadratic form on each
fiber, then the conic bundle is equivalent to an orthogonal bundle
(see \cite{R}). In this case there is already a definition of stability,
and we check in section \ref{orthogonalbundles} that it is a particular
case of our definition.

\begin{definition}
Consider a \tube\ $(\CrE,Q)$ and a subbundle $\CrE'$ of $\CrE$ of rank
$r'$. Let $x$ be a general point in $X$. If $\CrF_1$ and $\CrF_2$ are 
subbundles of $\CrE$, we denote by $\CrF_1 \CrF_2$ the subbundle of 
$\Sym^2 \CrE$
generated by elements of the form $f_1 f_2$ where $f_1$ and $f_2$ 
are local sections of $\CrF_1$ and $\CrF_2$. We define a function
$c_Q(\CrE')$ as follows:
$$
c_Q(\CrE')=\left\{
\begin{array}{lcl}
2, & \text{if} & Q|_{\CrE'\CrE'}\neq 0\\
1, & \text{if} & Q|_{\CrE'\CrE}\neq 0 = Q|_{\CrE'\CrE'} \\
0, & \text{if} & Q|_{\CrE'\CrE}=0.
\end{array}
\right.
$$
\end{definition}

Sometimes it will be convenient to write this type of conditions on
$Q$ in matrix form. Choosing a basis compatible with the filtration
$\CrE' \subset \CrE$ these three cases can be expressed as follows
$$
\left(
\begin{array}{cc}
\times & \cdot \\
\cdot  & \cdot 
\end{array}
\right)
,\qquad
\left(
\begin{array}{cc}
0      & \times \\
\times & \cdot 
\end{array}
\right)
,\qquad
\left(
\begin{array}{cc}
0   & 0 \\
0   & \times 
\end{array}
\right)
,\qquad
$$
where $\times$ means that that block is nonzero, $0$ means that it is
zero and $\cdot$ means that it can be anything.

\begin{definition}
\label{gcritical}
Let $(\CrE,Q)$ be a \tube. We say that two subbundles $\CrE_1 \subset \CrE_2
\subset \CrE$ give a critical filtration of $(\CrE,Q)$, if $\rk(\CrE_1)=1$,
$\rk(\CrE_2)=2$, $\rk(\CrE)=3$, $Q|_{\CrE_1 \CrE_2}=0$, 
and $Q|_{\CrE_1 \CrE}\neq 0 
\neq Q|_{\CrE_2 \CrE_2}$.
\end{definition}

The fact that $\CrE_1 \subset \CrE_2 \subset \CrE$ is a critical filtration of
$(\CrE,Q)$ means that for a generic point $x\in X$, the conic $Q_x$
defined by $Q$ on the fibre of $\PP(\CrE)$ over $x$ is smooth, the point
defined by $\CrE_1$ is in the conic and the line defined by $\CrE_2$ is
tangent to the conic. In matrix form with a basis adapted to the
filtration $\CrE_1\subset \CrE_2 \subset \CrE$ this can be expressed
as
$$
Q=\left(
\begin{array}{ccc}
0 & 0 & \times \\
0 & \times & \cdot \\
\times & \cdot & \cdot
\end{array}
\right)
$$
Later on (definition \ref{gcritfiltvectspac}) 
we will introduce a similar definition for filtrations of vector
spaces.

Now we are ready to define the notion of stability. We will only
define it for $r\leq 3$. As it is usual when one is
working with vector bundles with extra structure, this notion will
depend on a positive rational number $\tau$. We could as well take
$\tau$ to be a real number, but this wouldn't give anything new
because when we vary $\tau$ the stability of a \tube\ can only
change at rational values of $\tau$.  

We follow the notation of \cite{H-L}: Whenever the word '(semi)stable'
appears in a statement with the symbol '$(\leq)$', two statements should
be read. The first with the word 'stable' and strict inequality, and the
second with the word 'semistable' and the relation '$\leq$'.

\begin{definition}

Let $\tau$ be a positive rational number. Let $(\CrE,Q)$ be a \tube\
with $r\leq 3$. We say that $(\CrE,Q)$ is (semi)stable with respect to
$\tau$ if the following conditions hold

(ss.1) If $\CrE'$ is a proper subbundle of $\CrE$, then
$$
\frac{\deg(\CrE')-c_Q(\CrE')\tau}{\rk(\CrE')} (\leq) \frac{\deg(\CrE) -
2\tau }{r}.
$$

(ss.2) If $\CrE_1\subset \CrE_2 \subset \CrE$ is a critical filtration, then
$$
\deg(\CrE_1)+\deg(\CrE_2) (\leq) \deg(\CrE).
$$
\end{definition}
 
Note that condition (ss.1) is reminiscent of the stability conditions
for vector bundles with extra structure in the literature, but
condition (ss.2) is new. It is due to the fact that in a \tube, $Q$ is
a nonlinear object. So far all objects that have been considered were
linear, and this is why this kind of conditions didn't appear. This
nonlinearity is responsible for the fact that the proof is more
involved, and we have to consider only \tubes\ with $r\leq 3$.
For higher $r$ we expect to have more conditions of the form (ss.2).

\begin{lemma}
\label{gautom}
Let $(\CrE_1,Q_1)$ and $(\CrE_2,Q_2)$ be stable \tubes\ of the same type
$(r,d,\CrL)$. Then any nontrivial morphism $\phi:(\CrE_1,Q_1) \to
(\CrE_2,Q_2)$ is an isomorphism, and furthermore it is a scalar 
multiple of identity.
\end{lemma}

\begin{proof}
Assume that $\phi$ is nontrivial. Let $f:\CrE_1 \to \CrE_2$ be the
corresponding morphism of sheaves. Consider the subsheaves $\CrE'=\ker f$
of $\CrE_1$ and $\CrE''=\im f$ of $\CrE_2$. Assume $\CrE' \neq 0$. By
commutativity of the diagramme
$$
\CD
\Sym^2 \CrE' @>>>   \Sym^2 \CrE_1 @>\Sym^2 f>> \Sym^2 \CrE_2 \\
 @V{Q_1}VV        @V{Q_1}VV             @V{Q_2}VV \\
    \CrL       @=      \CrL        @>g>>          \CrL
\endCD
$$
we have that $c_{Q_1}(\CrE')=0$, and then by stability
$$
\frac{\deg(\CrE')}{\rk(\CrE')} < \frac{d-2 \tau}{r} <\frac{\deg(\CrE'')-2\tau}
{\rk(\CrE'')} \leq \frac{\deg(\CrE'')-c_{Q''}(\CrE'')\tau}{\rk(\CrE'')}
<\frac{d-2\tau}{r}
$$
which is a contradiction. Then $\CrE'=0$ and $f$ is an isomorphism.
Now let $x\in X$ be a point, and let $\lambda$ be an eigenvalue of $f$
at the fibre over $x$. Then $h=f-\lambda \id_{\CrE_1}$ is not surjective
at $x$, hence $h$ cannot be an isomorphism and then $h=0$.

\end{proof}

A flat family of (semi)stable \tubes\ of type $(r,d,\CrL)$ parametrized by a 
scheme $T$ is a triple $(\CrE_T,Q_T,\CrN)$ where $\CrE$ is a
vector bundle on $X\times T$, flat over $T$, that restricts to a 
vector bundle of rank
$r$ and degree $d$ on each fibre $X\times t$, and $Q_T$ is a 
morphims $Q_T:\Sym^2 \CrE_T \to p^*_X \CrL \otimes p^*_T \CrN$ 
where $\CrN$ is a line
bundle on $T$, and this morphim restricts to (semi)stable \tubes\ on
each fibre. Two families $(\CrE_T,Q_T,\CrN)$ and
$(\CrE'_T,Q'_T,\CrN')$ will be considered equivalent if there is a line
bundle $\CrM$ on $T$, an isomorphism $f:\CrE_T\otimes p^*_T \CrM \to
\CrE'_T$ and a commutative diagramme
$$
\CD
\Sym^2 \CrE_T\otimes p^*_T \CrM^2  @>\Sym^2 f>> \Sym^2 \CrE'_T \\
@V{Q}VV             @V{Q'}VV \\
p^*_X\CrL\otimes p^*_T\CrN \otimes p^*_T\CrM^2 @>{\isom}>> 
p^*_X\CrL\otimes p^*_T\CrN'
\endCD
$$
Let $\FM_\tau (r,d,\CrL)^\natural$ (resp. $\overline \FM_\tau (r,d,\CrL)
^\natural$) be
the functor that sends a scheme $T$ to the set of flat families of 
stable (resp. semistable) \tubes\
of type $(r,d,\CrL)$ parametrized by $T$. The moduli space for this
functor will be denoted by $\FM_\tau (r,d,\CrL)$ (resp. $\overline 
\FM_\tau (r,d,\CrL)$).

\smallskip
\noindent
\textbf{Theorem I.}
\label{gmaintheorem}
\textit{Let $X$ be a Riemann surface. Let $\tau>0$ be a rational number. There
exist a projective coarse moduli space $\overline \FM_\tau (r,d,\CrL)$ of
semistable \tubes\ with respect to $\tau$ of fixed type $(r,d,\CrL)$. The
closed points of $\overline \FM_\tau (r,d,\CrL)$ correspond to
S-equivalence classes of \tubes . There is an open set $\FM_\tau
(r,d,\CrL)$ corresponding to stable \tubes. This open set is a fine moduli
space of stable \tubes. Points in this open set
correspond to isomorphism classes of \tubes.}

For a definition of S-equivalence, see subsection \ref{gsequiv}

At the same time we wrote this article, I. Mundet i Riera found the
conditions for existence of solutions to a generalization of the 
vortex equation associated to Kaehler fibrations. 
As expected, the condition he finds is, 
in the case of conic bundles, the same as the condition we have found for 
stability. This is explained in the appendix.

\section{GIT Construction}
\label{ggit}

In this section we will construct the moduli space of semistable
\tubes. This construction is based on the ideas of Simpson
for the construction of the moduli space of semistable sheaves 
(\cite{S}). We will follow closely the paper \cite{K-N} of King 
and Newstead and the paper  \cite{H-L} of Huybrechts
and Lehn. In \ref{gboundnesstheorems} we prove some boundness
theorems that are needed later, and in \ref{gconstruction} we give
the construction of the moduli space and prove the semistability
condition. The base field $k$ can be any algebraically closed field of
characteristic zero, but we are mainly interested in $\CC$.

\subsection{Boundness theorems}
\label{gboundnesstheorems}

\begin{proposition}
\label{gboundnessb}
Let $X$ be a genus $g$ curve. Let $\calS$ be a set of vector bundles
on $X$ with degree $d$ and rank $r$. Assume that there is a constant
$b$ such that if $\CrE \in \calS$ and $\CrE'$ is a nonzero subsheaf of $\CrE$,
then 
$$
\mu(\CrE')=\frac{\deg (\CrE')}{\rk (\CrE')} \leq b.
$$
Then there is a constant $m_0$ such that if $m\geq m_0$, for all $\CrE
\in \calS$, we have $h^1(\CrE(m))=0$ and $\CrE(m)$ is generated by global
sections. Hence $\calS$ is bounded.
\end{proposition}

\begin{proof}
Let $x$ be a point of the curve $X$ and $\CrE\in \calS$. The exact sequence
$$
0 \to \CrE(m)\otimes \SO_X(-x) \to \CrE(m) \to \CrE(m)|_x \to 0
$$
gives that if $h^1(\CrE(m)\otimes \SO_X(-x))=0$ for all $x\in X$, then 
$\CrE(m)$ is
generated by global sections and $h^1(\CrE(m))=0$.

Assume that $h^1(\CrE(m)\otimes \SO_X(-x)) \neq 0$. Then by Serre 
duality there is a
nonzero morphism $\CrE(m)\otimes \SO_X(-x) \to \CrK_X$, where $\CrK_X$
is the canonical divisor. This gives an 
effective divisor
$D$ on $X$ and an exact sequence
$$
0 \to \CrE'(m) \to \CrE(m) \to \CrK_X(x-D) \to 0.
$$
Let $d'=\deg(\CrE')$. We have $\rk (\CrE')=r-1$. Then
$$
d'=(1-r)m+d+rm-(2g-1-\deg (D))\geq d-2g+1+m
$$
On the other hand, by hypothesis $d'\leq (r-1)b$, and combining both
inequalities we get 
$$
m\leq (r-1)b-d+2g-1.
$$
Then if we take $m_0>(r-1)b-d+2g-1$, for any $m\geq m_0$ and $x\in X$
we will have $h^1(\CrE(m)\otimes \SO_X(-x))=0$, thus $\CrE(m)$ is 
generated by global
sections and $h^1(\CrE(m))=0$. By standard methods using the Quot
scheme, this implies that $\calS$ is bounded.

\end{proof}

\begin{corollary}
\label{gboundcor1}
The same conclusion is true for the set of vector bundles
$\CrE$ occurring in semistable \tubes\ $(\CrE,Q)$ of fixed type.
The constant $m_0$ depends on $X$, $\tau$, $r$ and $d$, but not on
$\CrL$.
\end{corollary}

\begin{proof}
By condition (ss.1) we have that for every subsheaf $\CrE'$ of $\CrE$
$$
\frac{\deg (\CrE')}{\rk (\CrE')} \leq \frac{d-2\tau}{r}+\frac{c_Q
(\CrE') \tau}{\rk (\CrE')} \leq \frac{d-2\tau}{r}+2\tau.
$$
Take $b=\frac{d-2\tau}{r}+2\tau$ and apply the proposition.

\end{proof}

\begin{corollary}
\label{gboundcor2}
Let $\calS$ be the set of semistabilizing sheaves, i.e. sheaves $\CrE'$,
$\CrE_1$, $\CrE_2$ that give equality in condition (ss.1) or (ss.2). Then
the conclusions of proposition \ref{gboundnessb} are also true for
$\calS$.
\end{corollary}

\begin{proof}
By semistability, the slope of a subsheaf of a sheaf in $\calS$ is
bounded. On the other hand there are only a finite number of
possibilities for rhe rank and degree of a sheaf in $\calS$, then we
can apply proposition \ref{gboundnessb}.

\end{proof}

Now we will state two lemmas of King and Newstead
(\cite[lemma 2.2]{K-N} and \cite[corollary 2.6.2]{K-N}).

\begin{lemma}
Let $\CrE$ be a torsion free sheaf such that for all subsheaf $\CrF$ of $\CrE$,
$\mu(\CrF)\leq b$. If $b<0$, then $h^0(\CrE)=0$. If $b\geq 0$ then
$h^0(\CrE)\leq \rk(\CrE) (b+1)$.
\end{lemma}
\hfill $\Box$

\begin{lemma}
\label{gbound2}
Fix $R$, $b$, $k$. Then there exists an $n_0$ such that if $\calS$ 
is a set of torsion free sheaves with 

(i) $\rk(\CrE) \leq R$

(ii) $\mu(\CrF)\leq b$ for all nonzero subsheaves $\CrF$ of $\CrE$

(iii) For some $n\geq n_0$
$$
h^0(\CrE(n))\geq \rk(\CrE) (\chi(\SO_X (n))+k)
$$
Then the set $\calS$ is bounded.
\end{lemma}
\hfill $\Box$

\subsection{Construction and proof of main theorem}
\label{gconstruction}
\hfil
\medskip

Now we will give the GIT construction of the moduli space. We will
assign a point in a projective scheme $Z$ to a
\tube\ $(\CrE,Q)$ of fixed type $(r,d,\CrL)$. Let $P$ be the Hilbert
polynomial of $\CrE$, i.e. $P(m)=rm+d+r(1-g)$. We will assume that $m$ is
large enough so that corollaries \ref{gboundcor1} and \ref{gboundcor2} are
satisfied. Let $V$ be a vector space of dimension $p=P(m)$. Let $\SH$
be the Hilbert scheme $\Hilb(V\otimes \SO_X(-m),P)$ parametrizing
quotients of $V\otimes \SO_X(-m)$ with Hilbert polynomial $P$.
Let $l>m$ be an integer, $W=H^0(\SO_X (l-m))$, and $G$ be the 
Grassmannian $\Grass (V\otimes W,p)$ of quotients of $V\otimes W$ of
dimension $p$. For $l$ large enough we have embeddings
$$
\SH \to G \to \PP(\Lambda ^{P(l)}(V\otimes W))
$$
Let $B=H^0(\CrL)$ and $\SP=\PP(\Sym ^2(V^\vee \otimes B))$. Given $(\CrE,Q)$
and an isomorphism $V\isom H^0(\CrE(m))$ we get a point $\z$ in 
$\SH \times \SP$ as follows:

The vector bundle $\CrE(m)$ is generated by global sections (corollary
\ref{gboundcor1}), then we have a quotient
$$
q:V\otimes \SO_X(-m) \isom H^0(\CrE(m))\otimes \SO_X(-m) \surj \CrE.
$$
Denote by $\tilde q$ the point in $\SH$ corresponding to this
quotient. On the other hand, we get a point $\tilde Q$ in $\SP$ by the
composition
$$
\Sym ^2 V \isom \Sym^2 H^0(\CrE(m)) \to H^0(\Sym ^2\CrE(m)) \to H^0(\CrL(2m))=B
$$
Let $\tilde Q$ be a point in $\SP$. We will denote by $Q'$ a
representative of $\tilde Q$, i.e. $Q':\Sym^2 V \to B$. This gives 
(up to multiplication by a scalar) an evaluation
$$
\operatorname{ev}:\Sym^2 V \otimes \SO_X(-2m) \to B\otimes \SO_X(-2m)
\to \CrL.
$$
Let $Z$ be the closed subset of $\SH \times \SP$ of points $\z$ such
that (some multiple) of this evaluation map factors
through $\Sym^2 \CrE$
$$
\Sym^2 V \otimes \SO_X(-2m) \to \Sym^2 \CrE \to \CrL.
$$
The group $\slv$ acts in a natural way on $\SH \times \SP$.
A point in $Z$ will be called ``good'' if the quotient
$$
q:V\otimes \SO_X(-m) \surj \CrE
$$
induces an isomorphism $V \stackrel{\isom}\to H^0(\CrE(m))$, and $\CrE$ is 
torsion free. Note that a \tube\ $(\CrE,Q)$ gives a ``good'' point in $Z$
and conversely we can recover the \tube\ from the point, and two
``good'' points correspond to the same \tube\ iff they are in the same
orbit of the action of $\slv$. This action on $\SH \times \SP$ 
preserves the subscheme $Z$ and the subset of ``good'' points.

Let $\CrM$ be the line bundle on $\SH$ given by the embedding $\SH \to
\PP(\Lambda ^{P(l)}(V\otimes W))$. Embedd $Z$ in projective space with
$\SO_Z(n_1,n_2)=p_{\SH}^*\CrM ^{\otimes n_1} \otimes p_{\SP}^*
\SO_{\SP} (n_2)$
$$
Z \inj \PP(\Sym ^{n_1}[\Lambda ^{P(l)}(V\otimes W)] \otimes 
\Sym ^{n_2}[\Sym^2 V^\vee \otimes B])
$$
The group $\slv$ acts naturally on $\Sym ^{n_1}[\Lambda ^{P(l)}
(V\otimes W)] \otimes \Sym ^{n_2}[\Sym^2 V^\vee \otimes B]$, and
this gives a linearization for the action of $\slv$ on $Z$.

Now we will characterize the (semi)stable points of $Z$ under the
action of $\slv$ with the linearization induced by
$\SO_X(n_1,n_2)$. We will take
$$
\frac{n_2}{n_1}=\frac{P(l)-P(m)}{P(m)-2\tau} \tau.
$$

\textbf{Notation.} Given a point $\z\in Z$ and a subspace $V' \subset
V$ we denote by $\CrE_{V'}$ the image of $V'\otimes \SO_X(-m)$ under the
quotient $q:V\otimes \SO_X(-m) \to \CrE$. Note that $V' \subset
H^0(\CrE_{V'}(m))$, but in general they are not equal. If $\CrE' \subset \CrE$
is a subsheaf of $\CrE$ we have
$\CrE_{H^0(\CrE'(m))} \subset \CrE'$, with equality if $\CrE'(m)$ is generated by
global sections. Given a sheaf $\CrF$, we will denote
by $P_\CrF$ its Hilbert polynomial.

The following definition is analogous to definition \ref{gcritical}.

\begin{definition}
\label{gcritfiltvectspac}
Let $\z$ be a point in $Z$. Let
$V_1\subset V_2 \subset V_3=V$ be a filtration of $V$. Let $Q'_{ab}$ 
be the restriction of $Q':\Sym ^2 V \to B$ to
$V_a\otimes V_b$. We say that $V_1$, $V_2$ give a critical filtration of
$\z$, if $\rk(\CrE_{V_1})=1$, $\rk(\CrE_{V_2})=2$,
$Q'_{12}=0$, and $Q'_{13}\neq 0 \neq
Q'_{22}$.
\end{definition}

\begin{proposition}
\label{gmainprop}
For $l$ large enough the point $\z \in Z$ is
(semi)stable by the action of $\slv$ with respect to the linearization
by $\SO_Z(n_1,n_2)$ iff:

(*.1) If $V' \varsubsetneq V$ is a subspace of $V$, then
$$
\dim V' (n_1 P(l)+2n_2) (\leq) \dim V (n_1 P_{\CrE_{V'}}(l) +c_Q
(\CrE_{V'})n_2).
$$

(*.2) If $V_1\subset V_2 \subset V$ is a critical
filtration, then
$$
(\dim V_1 + \dim V_2)(n_1 P(l)+2 n_2) (\leq) \dim V
(n_1(P_{\CrE_{V_1}}(l)+ P_{\CrE_{V_2}}(l)) +2 n_2).
$$
\end{proposition}

\begin{proof}
We will apply the Hilbert-Mumford criterion: a point $(\tilde q,
\tilde Q)$ is (semi)stable iff for all one-parameter subgroup (1-PS)
$\lambda$ of $\slv$ we
have $\mu(\z,\lambda) (\leq) 0$, where
$\mu(\z,\lambda)$ is the minimum weight of the action
of $\lambda$ on $\z$.

Let $p=P(m)$. A 1-PS $\lambda$ of $\slv$ is equivalent to a basis
$\{v_1,\dots,v_p\}$ of $V$ and a weight vector
$(\gamma_1,\dots,\gamma_p)$ with $\gamma_i \in \ZZ$, $\gamma_1\leq
\dots \leq \gamma_p$, and $\sum \gamma_i=0$. The set 
$\SC$ of all weight vectors is a cone in
$\ZZ^p$. If a basis of $V$ has been chosen, then by a slight abuse of
notation we will denote $\mu(\z,\lambda)$ by $\mu(\z,\gamma)$, where
$\gamma \in \SC$.

We will choose a set of one-parameter subgroups, calculate
$\mu(\z,\lambda)$, and then imposing $\mu(\z,\lambda)(\leq)0$ we will
obtain necessary conditions for $\z$ to be (semi)stable.

Then we will show that the chosen set of one-parameter subgroups is
sufficient, in the sense that if we check that
$\mu(\z,\lambda)(\leq)0$ for all one-parameter subgroups in this set,
then the same will hold for any arbitrary one-parameter subgroup 
in $\SC$.

We have $\mu(\z,\lambda)=n_1 \mu(\tq,\lambda)+n_2 \mu(\tQ,\lambda)$, 
where $\mu(\tq,\lambda)$ (resp. $\mu(\tQ,\lambda)$) is the minimum
weight of the action of $\lambda$ on $\tq\in \SH$ (resp. $\tQ\in \SP$).
Fix a basis $\{v_1,\dots,v_p\}$ of $V$. Define $\varphi(i)=\dim q'(
\langle v_1,\dots,v_i \rangle \otimes W)$, where $q':V\otimes W 
\surj k^{P(l)}$ is
the quotient corresponding to the point $\tq \in \SH$. We have (see
\cite{S})
$$
\mu(\tq,\gamma)=\sum_{i=1}^p \gamma_i (\varphi(i)-\varphi(i-1)).
$$
On the other hand
$$
\mu(\tQ,\gamma)=\min _{i,j\in \{1,\dots,p\}} \{\gamma_i+\gamma_j :
Q'(v_i,v_j)\neq 0 \}.
$$
Note that $\mu(\tq,\gamma)$ is linear on $\gamma \in \SC$, but
$\mu(\tQ,\gamma)$ is not.

\bigskip
\noindent\textit{\textbf{GIT (semi)stable implies conditions (*)}}
\bigskip

Let $\z$ be a (semi)stable point in $Z$. Let $\{v_1,\dots,v_p\}$ be a 
basis of $V$. Define
$$
i_k=\min \{i: \rk(\CrE_{\langle v_1,\dots,v_i\rangle}) \geq k\}.
$$
Note that if $\z$ is ``good'', then the map $V\to H^0(\CrE(m))$ is an
isomorphism (in particular injective), and then $i_1=1$. Later on we
will see that for sufficiently large $m$, a semistable point is
``good'', but now we won't assume that $\z$ is ``good''. Define a
filtration of $V$
$$
V_1=\langle v_1,\dots,v_{i_1}\rangle \subset V_2=\langle
v_1,\dots,v_{i_2} \rangle 
\subset V_3=V.
$$
Let $Q'_{ab}$ be the restriction of $Q':\Sym ^2 V \to B$ to
$V_a\otimes V_b$. To calculate $\mu(\tQ,\lambda)$ we distinguish seven
cases.
$$
\begin{array}{rll}
1)& Q'_{11}\neq 0 & \mu(\tQ,\lambda)=2\gamma_{i_1} \\
2)& Q'_{11}=0,\; Q'_{12}\neq 0 &
\mu(\tQ,\lambda)=\gamma_{i_1}+\gamma_{i_2}  \\
3)& Q'_{12}=0,\; Q'_{13}\neq 0\neq Q'_{22} & \mu(\tQ,\lambda)=
\min(2\gamma_{i_2},\gamma_{i_1}+\gamma_{i_3}) \\
4)& Q'_{13}=0,\; Q'_{22}\neq 0 & \mu(\tQ,\lambda)= 2\gamma_{i_2} \\
5)& Q'_{22}=0,\; Q'_{13}\neq 0 & \mu(\tQ,\lambda)=
\gamma_{i_1}+\gamma_{i_3} \\
6)& Q'_{13}=Q'_{22}=0,\; Q'_{23}\neq 0 & \mu(\tQ,\lambda)=
\gamma_{i_2}+\gamma_{i_3} \\
7)& Q'_{23}=0,\; Q'_{33}\neq 0 & \mu(\tQ,\lambda)= 2\gamma_{i_3} \\
\end{array}
$$
Note that in all cases, except case 3, $\mu(\tQ,\lambda)$ is a linear
function of $\gamma \in \SC$.

First we will consider weight vectors of the form
\begin{eqnarray}
\gamma^{(i)}=(\overbrace{i-p,\dots,i-p}^i,\overbrace{i,\dots,i}^{p-i})
\; \; \;
(1\leq i < p).
\label{ggenerator1}
\end{eqnarray}
and define $V'=\langle v_1,\dots,v_i\rangle $ (it is clear that any 
subspace of $V$ can
be written in this form, after choosing an appropriate bases for $V$).

We have $\mu(\tq,\gamma^{(i)})=-p \varphi(i)+i \varphi(p)$. To obtain
a formula for $\mu(\tQ,\gamma^{(i)})$ we have to analyze each of the
seven cases. We will only work out the details for cases 2 and 3, the
remaining cases being similar to case 2.

In case 2 we have
$\mu(\tQ,\gamma^{(i)})=\gamma^{(i)}_{i_1}+\gamma^{(i)}_{i_2}$. Then,
according to the value of $i$ we have
$$
\mu(\tQ,\gamma^{(i)})=
\left\{
\begin{array}{lll}
2i &,i<i_1 &(\rk(\CrE_{V'})=0)\\
2i-p &,i_1\leq i < i_2&(\rk(\CrE_{V'})=1)\\
2i-2p &,i_2 \leq i&(\rk(\CrE_{V'})\geq 2)\\
\end{array}
\right .
$$
In case 3 we have $\mu(\tQ,\gamma^{(i)})= \min (2\gamma^{(i)},
\gamma^{(i)}_{i_1}+\gamma^{(i)}_{i_3})$, hence
$$
\mu(\tQ,\gamma^{(i)})=
\left\{
\begin{array}{lll}
2i &,i<i_1 &(\rk(\CrE_{V'})=0)\\
2i-p &,i_1\leq i < i_2&(\rk(\CrE_{V'})=1)\\
2i-2p &,i_2 \leq i&(\rk(\CrE_{V'})\geq 2)\\
\end{array}
\right .
$$
Doing the calculation for the seven cases we check that in every case
we have 
$$
\mu(\tQ,\gamma^{(i)})=2i-c_Q(\CrE_{V_i})p
$$
Then $\mu(\z,\gamma^{(i)})(\leq)0$ gives 
$$
n_1(-p \varphi(i) +i \varphi(p))+ n_2(2i-c_Q(\CrE_{V'}))(\leq) 0
$$
If we vary $V'$ (allowing $V'=V$), the submodules $\CrE_{V'}$ are
bounded, so we can take $l$ large enough such that
$\varphi(i)=P_{\CrE_{V'}}(l)$. We have $i=\dim V'$, and $\varphi(p)=P(l)$,
and then we obtain condition (*.1).

To obtain condition (*.2), assume that we have subspaces $V_1 \subset
V_2 \subset V$ giving a critical filtration. Let $i=\dim V_1$ and
$j=\dim V_2 - \dim V_1$. Take a bases $\{v_1,\dots,v_p\}$ of $V$ adapted
to this filtration, i.e. such that $V_1=\langle v_1,\dots,v_i\rangle $ and
$V_2=\langle v_1,\dots v_{i+j}\rangle $. Consider the weight vector

\begin{eqnarray}
\lefteqn{\gamma^{(i)}+\gamma^{(i+j)}=} \nonumber  \\
& &(\overbrace{2i+j-2p,\dots,2i+j-2p}^{i},
\overbrace{2i+j-p,\dots,2i+j-p}^{j},
\overbrace{2i+j,\dots,2i+j}^{p-i-j}).\label{ggenerator2}
\end{eqnarray}
An easy computation then shows
$\mu(\tQ,\gamma^{(i)}+\gamma^{(i+j)})= 2(\dim V_1+\dim V_2 -\dim V)$.
On the other hand $\mu(\tq,
\gamma^{(i)}+\gamma^{(i+j)})= (\dim V_1+\dim V_2)P(l) -\dim
V(P_{V_1}(l) + P_{V_2}(l))$, and then
$\mu(\z,\gamma^{(i)}+\gamma^{(i+j)}) (\leq) 0$ gives condition (*.2).

\bigskip
\noindent\textit{\textbf{ Conditions (*) imply GIT (semi)stable}}
\bigskip

Now we have to show that the one-parameter subgroups that we have used
are sufficient. As we did before, we will fix an arbitrary base $V$,
and we consider the seven different cases. In all cases except 3,
$\mu(\tQ,\gamma)$, and hence $\mu(\z,\gamma)$, is a linear function of
$\gamma \in \SC$, and then to prove that $\mu(\z,\gamma)(\leq)0$ for
all $\gamma$ it is enough to check it on the generators $\gamma^{(i)}$
defined above (\ref{ggenerator1}).

In case 3 we have $\mu(\tQ,\gamma)= \min(2\gamma_{i_2} , \gamma_{i_1} +
\gamma_{i_3})$, hence it is no longer linear on $\gamma$, and it is
not enough to check the condition on the generators
$\gamma^{(i)}$. But it is a piecewise linear function. The cone $\SC$
of weights is divided in two cones
\begin{eqnarray*}
\SC^>=\{(\gamma_1,\dots,\gamma_p)\in \SC: 2\gamma_{i_2} \geq \gamma_{i_1} +
\gamma_{i_3}\} \\
\SC^<=\{(\gamma_1,\dots,\gamma_p)\in \SC: 2\gamma_{i_2} \leq \gamma_{i_1} +
\gamma_{i_3}\} 
\end{eqnarray*}
Observe that $\mu(\tQ,\gamma)$ is linear on each of these cones. We
will use the following lemma.

\begin{lemma}
\label{gconelemma}
Let $\SC$ be a cone in $\ZZ^p$, let $\gamma^{(i)}$ be a set of
generators of $\SC$, i.e. $\SC=(\oplus_i \QQ^+ \gamma^{(i)}) \cap
\ZZ^p$. Let $A:\ZZ^p \to \QQ$ be a linear function such that
$A(\gamma^{(i)}) \in \{1,0,-1\}$. Let $\SC^>$ be the subcone $\{v\in
\SC: A(v)\geq 0\}$. Then the set of vectors 
$$
v_{i,j}=
\left\{
\begin{array}{ll}
\gamma^{(i)} &,\; A(e_i)\geq 0 \\
\gamma^{(i)}+\gamma^{(i+j)} &,\; A(e_i)=-1,\; A(e_{i+j})=1 \\
0 &,\; \text{otherwise.}\\
\end{array}
\right .
$$
generate $\SC^>$.
\end{lemma}
\hfill $\Box$

We apply this lemma with $A(\gamma)=(2\gamma_{i_2}- \gamma_{i_1}-
\gamma_{i_3})/p$ (and then with the negative of this, for $\SC^<$), 
and we obtain a set of generators for $\SC^>$ and $\SC^<$. But all
these vectors are either of the form $\gamma^{(i)}$ with $1\leq i < p$, 
or of the form $\gamma^{(i)}+
\gamma^{(i+j)}$ with $\rk(\CrE_{\langle v_1,\dots,v_i\rangle })=1$ and
$\rk(\CrE_{\langle v_1,\dots,v_{i+j}\rangle })=2$, and we have already
considered
them.

\end{proof}

\begin{remark}
\textup{In the following propositions we will prove that conditions
(*) are equivalent to the stability conditions (s). Recall that $\dim
V=p=P(m)$, $\varphi(i)=P_{E_{V'}}(l)$, $\varphi(p)=P(l)$ and $\dim
V'=i$. The idea is to show that for $l\gg m \gg 0$, we can replace
$P(l)$ by $\rk(E)l$,  $P_{E_{V'}}(l)$ by $\rk(E')l$, $P(m)$ by
$\deg(E)+ rm$ and $\dim V'$ by $P_{E_{V'}}(m)$, and this by $\deg(E')+
rm$.}
\end{remark}

\begin{proposition}
\label{gtrans}
For $m$ and $l$ large enough we have that conditions (*) are equivalent
to:

(**.1) If  $V'\varsubsetneq V$ is a subspace of $V$, then
$$
r(\dim V' - c_Q(\CrE_{V'})) \leq \rk(\CrE_{V'})(\dim V -
2 \tau),
$$
and in case of equality we also require $\dim
V'(\leq) P_{\CrE_{V'}}(m)$.

(**.2) If $V_1\subset V_2 \subset V$ is a critical filtration, then
$$
\dim V_1 + \dim V_2 \leq \dim V
$$
and in case of equality we also require $\dim V_1 + \dim V_2 (\leq)
P_{\CrE_{V_1}}(m)+P_{\CrE_{V_2}}(m)$.
\end{proposition}

\begin{proof}
We rewrite (*.1) using 
$$
\frac{n_2}{n_1}=\frac{P(l)-P(m)}{P(m)-2\tau} \tau.
$$
We obtain
\begin{eqnarray*}
[(\dim V_1 - c_Q(\CrE_{V'})\tau)r - \rk(\CrE_{V'})(\dim V -
2\tau )](l-m) + \\
+(\dim V - 2 \tau) \dim V (\dim V' - P_{\CrE_{V'}})(m)
(\leq) 0.
\end{eqnarray*}
We have $(l-m)\gg 0$ and $m\gg 0$, hence $\dim V >2\tau$ and the 
result follows. Now we rewrite (*.2), using
$r=3$, $\rk(\CrE_{V_1})=1$ and $\rk(\CrE_{V_2})=2$.
\begin{eqnarray*}
3(\dim V_1 + \dim V_2 - \dim V)(l-m) + \\
+(\dim V -2\tau)\dim V (\dim
V_1 + \dim V_2 - P_{\CrE_{V_1}}(m)-P_{\CrE_{V_2}}(m))(\leq)0,
\end{eqnarray*}
and the result follows.

\end{proof}

\begin{proposition}
\label{giden}
For $m$ and $l$ large enough, we have

(i) If $(\CrE,Q)$ is a (semi)stable \tube, 
then the corresponding point $\z$ in $Z$ is $GIT$ (semi)stable under
the action of $\slv$.

(ii) If $\z\in Z$ is a $GIT$ semistable point, then $\tq$ is
``good'' and $h^1(\CrE(m))=0$.

(iii) If $\z\in Z$ is a $GIT$ (semi)stable point, then the
corresponding \tube\ $(\CrE,Q)$ is (semi)stable.

Note that thanks to (ii), in (iii) we know that $\CrE$ is torsion free.
\end{proposition}

\begin{proof}
We will proof the three items in three steps

\bigskip
\noindent\textit{\textbf{Step 1. (Semi)stable \tube\ $\Rightarrow$ 
GIT (semi)stable $\z$}}
\bigskip

We will use proposition \ref{gtrans}. We will start checking (**.1).
Let $\calS$ be the set of vector bundles $\CrE'$ that are subsheaves of
bundles $\CrE$ occurring in semistable \tubes. It satisfies hypothesis
(i) and (ii) of lemma \ref{gbound2} with $R=3$ and $b=\frac{d-2 \tau}
{r} +2 \tau$. Let $k=\frac{d-2\tau}{r}$, $n$ large enough, so that
propositions \ref{gbound2}, \ref{gboundcor1} and \ref{gboundcor2} hold,
and let $\calS_n$ be the subset of $\calS$ consisting of bundles $\CrE'$
that satisfy hypothesis (iii) of lemma \ref{gbound2}. Then the
set $\calS_n$ is bounded. Taking $m>n$ large enough we then have
$h^1(\CrE'(m))=0$ for $\CrE' \in \calS_n$. In other words,
\begin{eqnarray}
h^0(\CrE'(m))=\rk(\CrE')(\chim + \frac{\det (\CrE')}{\rk(\CrE')}),\;\;\;\;\;
\text{for}\;\; \CrE'\in \calS_n 
\end{eqnarray}
On the other hand, we still have
\begin{eqnarray}
\label{gineq1}
h^0(\CrE'(m))<\rk(\CrE')(\chim + \frac{d-2\tau}{r}),\;\;\;\;\;
\text{for}\;\; \CrE'\in \calS \setminus \calS_n
\end{eqnarray}
Let $V'$ be a subspace of $V$, and $\CrE_{V'}$ the corresponding
sheaf. If $\CrE_{V'}$ belongs to $\calS_n$, we get that condition (ss.1)
implies (**.1), because
\begin{eqnarray*}
\frac{\dim V' - c_Q(\CrE_{V'})\tau}{\rk(\CrE_{V'})} \leq 
\frac{h^0(\CrE_{V'}(m))-c_Q(\CrE_{V'})\tau}{\rk(\CrE_{V'})} = & \\
=\frac{\deg (\CrE_{V'})-c_Q(\CrE_{V'})\tau}{\rk(\CrE_{V'})}+\chim
(\leq) & \\
(\leq) \frac{d-2\tau}{r}+\chim = \frac{\dim V - 2\tau}{r}
\end{eqnarray*}
and $\dim V' \leq h^0(\CrE_{V'}(m))=P_{\CrE_{V'}}(m)$, because 
$h^1(\CrE'(m))=0$.
On the other hand, if $\CrE_{V'}$ belongs to $\calS\setminus\calS_n$,
inequality (\ref{gineq1}) implies (**.1)
$$
\frac{\dim V' - c_Q(\CrE_{V'})}{\rk(\CrE_{V'})} \leq 
\frac{h^0(\CrE_{V'}(m))}{\rk(\CrE_{V'})} <
\frac{d-2\tau}{r}+\chim=\frac{\dim V -2\tau}{r}
$$
In both cases,
if inequality (ss.1) is strict, then inequality (**.1) is
also strict. But assume that there is a semistabilizing subsheaf
$\CrE'$ of $\CrE$ (i.e. giving equality in (ss.1)). By corollary
\ref{gboundcor2}, $\CrE'(m)$ is generated by global sections. Let
$V'=H^0(\CrE'(m)) \subset H^0(\CrE(m))=V$. Then $\CrE'=\CrE_{V'}$, 
and we have
\begin{eqnarray*}
\lefteqn{\frac{\dim V' -c_Q(\CrE_{V'})\tau}{\rk(\CrE_{V'})}  = 
\frac{\deg(\CrE_{V'})-c_Q(\CrE_{V'})}{\rk(\CrE_{V'})}+\chim=} \\
& & \frac{d-2\tau}{r}+\chim=\frac{\dim V -2\tau}{r}
\end{eqnarray*}
and $\dim V'=P_{\CrE_{V'}}(m)$.

Now we will check condition (**.2). Let $\ST$ be the set of vector
bundles of the form $\CrE_1\oplus \CrE_2$ such that $\CrE_1 \subset \CrE_2 \subset
\CrE$ gives a critical filtration of a (semi)stable \tube\
$(\CrE,Q)$. Hypothesis (i) and (ii) of lemma 
\ref{gbound2} are satisfied with $R=3$ and
$b=\frac{d-2\tau}{r}+2\tau$. Let $k=d/3$, and $n$ large
enough. Let $\ST_n$ be the subset of $\ST$ consisting of vector
bundles $\CrE_1\oplus \CrE_2$ satisfying hypothesis (iii). Then $\ST_n$ is
bounded, and taking $m$ large enough we have $0=h^1((\CrE_1\oplus
\CrE_2)(m))=h^1(\CrE_1(m))+h^1(\CrE_2(m))$ for $\CrE_1\oplus \CrE_2 \in \ST_n$.
Hence for $\CrE_1\oplus \CrE_2\in \ST_n$,
\begin{eqnarray}
h^0(\CrE_1(m))+h^0(\CrE_2(m))=3\chim + \deg (\CrE_1)+\deg(\CrE_2)
\end{eqnarray}
On the other hand, for $\CrE_1\oplus \CrE_2\in \ST\setminus \ST_n$ we 
still have
\begin{eqnarray}
\label{gineq2}
h^0(\CrE_1(m))+h^0(\CrE_2(m))<3\chim + d
\end{eqnarray}
Let $V_1 \subset V_2 \subset V$ be a critical filtration of $V$. If
$\CrE_{V_1}\oplus \CrE_{V_2} \in \ST_n$, we get that (ss.2) implies (**.2),
because 
\begin{eqnarray*}
\dim V_1 + \dim V_2 \leq
h^0(\CrE_{V_1}(m))+h^0(\CrE_{V_2}(m))= \\
=3\chim + \deg (\CrE_{V_1})+ \deg(\CrE_{V_2}) (\leq) 
 3\chim + \deg(\CrE)=\dim V
\end{eqnarray*}
and also $\dim V_1 + \dim V_2 \leq h^0(\CrE_{V_1}(m)) + h^0(\CrE_{V_2}(m)) =
P_{\CrE_{V_1}}(m) + P_{\CrE_{V_2}}(m)$.

On the other hand, if $\CrE_{V_1}\oplus \CrE_{V_2} \in \ST\setminus \ST_n$,
inequality (\ref{gineq2}) implies (**.2)
$$
\dim V_1 + \dim V_2 \leq h^0(\CrE_{V_1}(m)) + h^0(\CrE_{V_2}(m)) < 3\chim +d
=\dim V
$$
In both cases, if inequality (ss.2) is strict, also (**.2) is
strict. But assume that we have subsheaves $\CrE_1 \subset \CrE_2 \subset \CrE$
giving a critical filtration of a semistable \tube\ $(\CrE,Q)$. By lemma
\ref{gboundcor2} $\CrE_1(m)$ and $\CrE_2(m)$ are generated by global sections
and $h^1(\CrE_1(m))=h^1(\CrE_2(m))=0$. Taking $V_1=H^0(\CrE_1(m))$ and $V_2=
H^0(\CrE_2(m))$ we have $\CrE_{V_1}=\CrE_1$ and $\CrE_{V_2}=\CrE_2$, hence
\begin{eqnarray*}
\dim V_1 +\dim V_2 = \deg(\CrE_1)+\deg(\CrE_2) +3\chim =\\
= \deg(\CrE)+3\chim =\dim V
\end{eqnarray*}
and also $\dim V_1 + \dim V_2 = P_{\CrE_{V_1}}(m)+P_{\CrE_{V_2}}(m)$

\bigskip
\noindent\textit{\textbf{Step 2. $\z$ GIT (semi)stable $\Rightarrow$ 
$h^1(\CrE(m))=0$ and $\tq$ good}}
\bigskip

If $h^1(\CrE(m))\neq 0$, then by Serre duality $\Hom (\CrE,\CrK_X)\neq 0$. Take
$\psi \in \Hom (\CrE,\CrK_X)$. The composition $V\otimes \SO_X \to \CrE(m) \to
\CrK_X$ gives a linear map
$$
f: V \to H^0(\CrK_X).
$$
Let $U$ be the kernel of $f$. We have $\dim U \geq \dim V - \dim
H^0(\CrK_X) = p-g$. Then by (semi)stability of $\z$ we have
$$
r(p-g-c_Q(\CrE_U)\tau) \leq r(\dim U- c_Q(\CrE_U)\tau) 
(\leq) \rk (E_U) ( p- 2\tau),
$$
$$
(r-\rk(\CrE_U))p \leq r(g+c_Q(\CrE_U)\tau) - \rk(\CrE_U) 2 \tau
$$
We have $\rk(\CrE_U)\leq r$. Then if $m$ is large enough the inequality
forces $r=\rk(\CrE_U)$. By definition of $U$ we have $\CrE_U(m)\subset
\ker \psi$, then $\rk(\ker \psi)=r$, $\rk(\im \psi)=0$, and then
$\psi=0$ because $\CrK_X$ is torsion free. We conclude that (for $m$
large enough) $h^1(\CrE(m))=0$.

Then $\dim V=p=h^0(\CrE(m))$, and to show that $\tq$ is ``good'' it is
enough to show that the induced linear map
$$
V \to H^0(\CrE(m))
$$
is injective. Let $V'$ be the kernel. Then we have
$\rk(\CrE_{V'})=0$. By semistability we have (**.1)
$$
r(\dim V' - c_Q(\CrE_{V'})) \leq 0,
$$
but $c_Q(\CrE_{V'})=0$, and then $\dim V'$ must be zero.

To show that $\CrE$ is torsion free, let $\CrT\subset \CrE$ be the torsion
subsheaf. We have $V \isom H^0(\CrE(m))$, and then $U=H^0(\CrT(m))$ is a
subspace of $V$. The associated sheaf $\CrE_U$ has rank equal to zero,
and arguing as above we get $U=0$.

\bigskip
\noindent\textit{\textbf{Step 3. GIT (Semi)stable $\z$ $\Rightarrow$ 
(semi)stable \tube}}
\bigskip

By the previous step we know that we can choose $m$ large enough so
that $\tq$ is ``good''. We will check first (ss.1). Let $\CrE'$ be a
subsheaf of $\CrE$. Define $V'=H^0(\CrE'(m))$. We have $\CrE_{V'}\subset \CrE'$,
$\rk (\CrE_{V'})\leq \rk(\CrE')$, $\dim V'\geq P_{\CrE'}(m)$, and
$c_Q(\CrE') \geq c_Q(\CrE_{V'})$. Then
\begin{eqnarray*}
r(P_{\CrE'}(m)-c_Q(\CrE')\tau ) \leq 
r(\dim V' - c_Q(\CrE_{V'}) \tau ) (\leq)\\
(\leq) \rk(\CrE_{V'})(\dim V -2 \tau)
\leq \rk(\CrE')(\dim V -2 \tau).
\end{eqnarray*}
Note that if (**.1) is strict, then also (ss.1) is strict. But assume
that there is a subspace $V'\subset V$ that is semistabilizing,
i.e. both conditions in (**.1) are equalities. Then
\begin{eqnarray*}
\frac{\deg (\CrE_{V'})-c_Q(\CrE_{V'}) \tau}{\rk(\CrE_{V'})} =
\frac{\dim V' -c_Q (\CrE_{V'}) \tau}{\rk(\CrE_{V'})} -\chim =\\
=\frac{\dim V - 2 \tau}{r} -\chim =
\frac{\deg(\CrE) - 2 \tau}{r},
\end{eqnarray*}
and we get that (ss.1) for $\CrE_{V'}$ also gives equality.

Now we are going to check (ss.2). As in step 1, consider the set $\ST$
of vector bundles of the form $\CrE_1 \oplus \CrE_2$ such that $\CrE_1 \subset
\CrE_2 \subset \CrE$ gives a critical filtration. We have already proved
condition (ss.1), thus hypothesis  (i) and (ii) of lemma \ref{gbound2}
are again satisfied. Then, as in step 1, we define the subset $\ST_n
\subset \ST$, and taking $m$ large enough we can assume that the vector
bundles $\CrE_1$ and $\CrE_2$ are generated by global sections if
$\CrE_1\oplus \CrE_2 \in \ST$.

Let $\CrE_1 \subset \CrE_2 \subset \CrE$ be a
critical filtration of $(\CrE,Q)$. Let $V_1=H^0(\CrE_1(m))$ and
$V_2=H^0(\CrE_2(m))$. 
If $\CrE_{V_1}\oplus \CrE_{V_2} \in \ST_n$, then $\CrE_{V_1}$ and $\CrE_{V_2}$ are
generated by global sections and then $\CrE_{V_1}=\CrE_1$, $\CrE_{V_2}=\CrE_2$,
and $V_1\subset V_2 \subset V$ is a critical
filtration of $V$ and (**.2) holds
\begin{eqnarray*}
\deg (\CrE_1) + \deg(\CrE_2) =
\dim V_1 + \dim V_2 -3\chim (\leq)\\
(\leq)
\dim V - 3\chim = \deg(\CrE).
\end{eqnarray*}
On the other hand, if $\CrE_{V_1}\oplus \CrE_{V_2} \in \ST \setminus \ST_n$,
inequality (\ref{gineq2}) implies (**.2)
$$
\deg (\CrE_1) + \deg(\CrE_2) \leq 
h^0(\CrE_1(m)) + h^0(\CrE_2(m)) -3\chim <\deg(\CrE).
$$

Note that if (**.2) is strict then (ss.2) is also
strict.
But assume that there is a semistabilizing critical sequence $V_1
\subset V_2 \subset V$, i.e. a critical sequence giving equality in
both conditions of (**.2). Then
\begin{eqnarray*}
\deg(\CrE_{V_1})+\deg(\CrE_{V_2}) = \dim V_1 + \dim V_2 -3\chim =\\
=\dim V -3 \chim = \deg \CrE,
\end{eqnarray*}
and we also get an equality in (ss.2).

\end{proof}

Once we have established proposition \ref{giden} and lemma
\ref{gautom} we can prove theorem I using standard techniques. We
follow closely \cite{H-L}.

\begin{proof2}\textit{ of theorem I.}
Let ${\overline \FM_\tau (r,d,\CrL)}$ (resp. $\FM_\tau (r,d,\CrL)$) be
the GIT quotient of $Z$ (resp. $Z^s$) by $\slv$.

First we construct a universal family on $Z^{ss}$ using the universal
families of the Quot scheme $\SQ$ and on $\SP=\PP(\Sym^2V^\vee \otimes
B)$ (we think of $\SP$ as the Grassmannian of one dimensional
subspaces of $\Sym^2 V^\vee \otimes B$, and hence the universal subbundle of
subspaces is $\SO_{\SP}(-1)$).

Recall that $Z^{ss}$ is in $\SH \times \SP$. The universal quotient
$\CrE_\SH$ on
$\SH \times X$ pulls back to a vector bundle $\CrE_{Z^{ss}}=
p^*_{\SH \times \SQ} \CrE_\SH$ on $Z^{ss}$. On the other hand 
the universal subbundle on
$\SP \times X$ gives a morphism $\Sym^2 V \otimes p^*_X \SO_X(-2m) 
\to p^*_X\CrL \otimes \SO_{\SP}(1)$ of sheaves over $Z^{ss} \times
X$. By the definition of $Z$, there
is a line bundle $\CrN$ on $Z^{ss}$ such that this last morphism
factors and gives $Q_{Z^{ss}}:\Sym^2\CrE_{Z^{ss}} \to p^*_X \CrL
\otimes p^*_{Z^{ss}}\CrN$. 
Note that the line bundle $\CrN$ is needed because the factorization
on $Z$ is only up to scalar multiplication.
The triple 
$(\CrE_{Z^{ss}},Q_{Z^{ss}},\CrN)$ is a universal \tube.

Given a family 
$(\CrE_T,Q_T)$ of \tubes\ parametrized by $T$, and using the universal
family on $Z^{ss}$, we obtain a morphism $T
\to \overline \FM_\tau (r,d,\CrL)$. This is done in the following way:
Let $m$ be large enough so that proposition \ref{giden} holds. Given a
family $(\CrE_T,Q_T,\CrN)$ of \tubes\ parametrized by $T$, consider the
locally free sheaf $\SV={p_T}_*(\CrE_T\otimes p^*_X\SO_X(m))$, and note
that $p^*_T\SV \otimes p^*_X\SO_X(-m) \to \CrE_T$ is a surjection. Cover
$T$ with open sets $U_i$ such that there are isomorphisms $\phi_i:
V\otimes \SO_{U_i} \to \SV$. Then we have quotients $q_i:V\otimes
p^*_X \SO_X(-m) \surj \CrE_{U_i}$ and families of subspaces $\SO_{U_i} \inj
\Sym^2V^\vee \otimes B \otimes \SO_{U_i}$, and these give maps $U_i \to
Z^{ss}$. On the intersections $U_i \cap U_j$ this maps in general will
differ by the action of $\slv$, then they combine to give a well
defined morphism $T \to \overline \FM_\tau (r,d,\CrL)$.

It is straightforward to check the universal property for
$\overline \FM_\tau (r,d,\CrL)$, and then $\overline \FM_\tau
(r,d,\CrL)$ is a coarse moduli space.

Now we will show that the universal family restricted to $Z^s$
descends to $\FM_\tau (r,d,\CrL)$, making it a fine moduli space. 
Applying Luna's
\'etale slice theorem \cite{L}, we can find an \'etale cover $U'$ of 
$\FM_\tau (r,d,\CrL)$ over
which there is a universal family $(\CrE'_{U'},Q'_{U'})$. Consider $U''=
U' \times _{\FM_\tau (r,d,\CrL)} U'$ and take an isomorphism $\Phi:
p^*_1(\CrE'_{U'},Q'_{U'}) \to p^*_2(\CrE'_{U'},Q'_{U'})$ with the condition 
$p^*_1 Q'_{U'} = p^*_2 Q'_{U'} \circ \Sym^2 \Phi$. This isomorphism
exists and is unique by lemma \ref{gautom}, and then it satisfies the
cocycle condition of descend theory \cite[Chap. VII]{Mr}, and hence the
family $(\CrE'_{U'},Q'_{U'})$ descends to $\FM_\tau (r,d,\CrL)$.

\end{proof2}

\subsection{S-equivalence}
\label{gsequiv}
\hfil
\medskip

Let $(\CrE,Q)$ and $(\CrE',Q')$ be two nonisomorphic \tubes. If they
are strictly semistable, it
could still happen that the corresponding points in the moduli space
$\ssmoduli$ coincide. In this
case we say that they are S-equivalent (note that this is not the usual 
definition. Usually one defines two bundles as S-equivalent if the 
graded objects of their
Jordan-H\"older filtrations coincide, and then proves that S-equivalence 
classes corresponds to points of the moduli space). In this section, given a
strictly semistable \tube\ $(\CrE,Q)$, we will show
how to obtain a canonical representative $(\CrE^S,Q^S)$ of its 
S-equivalent class. In other words, given two semistable \tubes\ 
$(\CrE,Q)$ and $(\CrE',Q')$, they will be S-equivalent iff $(\CrE^S,
Q^S)$ is isomorphic to $({\CrE'}^S,{Q'}^S)$.

Let $(\CrE,Q)$ be a strictly semistable \tube. Then there exists at least
one ``semistabilizing object'', i.e. there exists either a subbundle
$\CrE' \subset \CrE$ that gives equality on (ss.1)(and then we say
that $\CrE'$ is a semistabilizing object of type I), or there is a critical
filtration $\CrE_1 \subset \CrE_2 \subset \CrE$ giving equality on
(ss.2)(and then we say that the filtration is a semistabilizing object
of type II). Choose one semistabilizing object. We define a new \tube\
$(\CrE_0,Q_0)$ as follows (it will depend on which semistabilizing object
we choose):

In the first case (corresponding to (ss.1)), the vector bundle is
defined to be $\CrE_0=\CrE' \oplus \CrE/\CrE'$ (note that if $\CrE$ 
is semistable and
$\CrE'$ gives equality on (ss.1), then $\CrE/\CrE'$ is torsion free). To define
$Q_0$, let $v$ and $w$ be local sections of $\CrE_0$ on an open set $U$. We
distinguish three cases:
$$
\text{If}\;\; c_Q(\CrE')=
\left\{
\begin{array}{ll}
2 \text{, then}\;
Q_0(v,w)=&
\left\{
\begin{array}{ll}
Q(v,w)&v,w \in \CrE'(U)\\
0&\text{otherwise}\\
\end{array}
\right .\\
1 \text{, then}\;
Q_0(v,w)=&
\left\{
\begin{array}{ll}
Q(v,w)&v\in \CrE'(U) \; \text{or}\; w\in \CrE'(U)\\
0&\text{otherwise}\\
\end{array}
\right .\\
0 \text{, then}\;
Q_0(v,w)=&
Q(v,w)\\
\end{array}
\right .
$$
In matrix form this can be written as follows
$$
\text{If} \; Q=
\left(
\begin{array}{cc}
\times & \cdot \\
\cdot  & \cdot
\end{array}
\right)
,\quad\text{then}\; Q_0=
\left(
\begin{array}{cc}
\times & 0 \\
0      & 0
\end{array}
\right)
$$
$$
\text{If} \; Q=
\left(
\begin{array}{cc}
0      & \times \\
\times & \cdot
\end{array}
\right)
,\quad\text{then}\; Q_0=
\left(
\begin{array}{cc}
0      & \times \\
\times & 0
\end{array}
\right)
$$
$$
\text{If} \; Q=
\left(
\begin{array}{cc}
0 & 0 \\
0 & \times
\end{array}
\right)
,\quad\text{then}\; Q_0=
\left(
\begin{array}{cc}
0  & 0 \\
0  & \times
\end{array}
\right)
$$
It is easy to see that this is well defined. 

In the second case
(corresponding to (ss.2)) we define the vector bundle to be $\CrE_0=\CrE_1
\oplus \CrE_2/\CrE_1 \oplus \CrE/\CrE_2$. Again let $v$ and $w$ be
local sections
of $\CrE_0$ on an open set $U$. Then we set
$$
Q_0(v,w)=
\left\{
\begin{array}{ll}
Q(v,w)&v \; \text{and}\; w \in \CrE_2(U)\\
Q(v,w)& v \; \text{or}\; w \in \CrE_1(U)\\
0&\text{otherwise,}\\
\end{array}
\right . 
$$
and in matrix form
$$
Q=
\left(
\begin{array}{ccc}
0      & 0      & \times \\
0      & \times & \cdot  \\
\times & \cdot  & \cdot  
\end{array}
\right)
\quad  Q_0=
\left(
\begin{array}{ccc}
0      & 0      & \times \\
0      & \times & 0 \\
\times & 0      & 0
\end{array}
\right)
.
$$
Again it is easy to see that this is well defined.

\begin{proposition}
\label{sequivcrit}
The \tube\ $(\CrE_0,Q_0)$ is also semistable. Furthermore, if we repeat
this process, eventually we will get a \tube\ that we will call
$(\CrE^S,Q^S)$ with the following properties

(i) $(\CrE^S,Q^S)$ is semistable, and if we apply this process to it with any
object we obtain an isomorphic \tube\ (i.e. this process stops).

(ii) $(\CrE^S,Q^S)$ only depends on the isomorphism class of $(\CrE,Q)$.

(iii) Two \tubes\ $(\CrE,Q)$ and $(\CrE',Q')$ are S-equivalent if and only if
$(\CrE^S,Q^S)$ is isomorphic to $({\CrE'}^S,{Q'}^S)$.
\end{proposition}

\begin{remark}
\textup{
The \tube\ $(\CrE^S,Q^S)$ is the analogue of the graded object $\gr(\CrE)$ 
of the Jordan-H\"older filtration of a semistable torsion-free sheaf.
Note that $\gr(\CrE)$ can also be obtained by a process similar to this.}
\end{remark}

\begin{proof}
We start with a general observation about GIT quotients. Let $Z$ be a
projective variety with a linearized action by a group $G$.
Two points in the open subset $Z^{ss}$ of semistable points are
S-equivalent (they are mapped to the same point in the moduli space)
if there is a common closed orbit in the closures (in $Z^{ss}$) of
their orbits. Let $z\in Z^{ss}$. Let $B$ be the unique closed orbit in
the closure $\overline{G\cdot z}$ in $Z^{ss}$ of its orbit $G \cdot
z$. Assume that $z$
is not in $B$. Then there exists a one-parameter 
subgroup $\lambda$ such that the limit $z_0=\lim_{t\to 0}
\lambda (t) \cdot z$ is in $\overline{G
\cdot z}\setminus G \cdot z$. Note that we must have
$\mu(z,\lambda)=0$ (otherwise $z_0$ would be unstable). 
Note that $G \cdot z_0 \subset \overline{G
\cdot z}\setminus G \cdot z$, and then $\dim G \cdot z_0 < \dim G
\cdot z$. Repeating this process with $z_0$ we then get a sequence 
of points that eventually stops and gives $\tilde z \in B$. Two points
$z_1$ and $z_2$ will then be S-equivalent if and only if after
applying this procedure to both of them the orbits of
$\tilde z_1$ and $\tilde z_2$ are the same.

We will use the notation introduced in subsection
\ref{gconstruction}. 
We will prove the proposition using the previous observation. The fact
that choosing a ``semistabilizing object'' of $(\CrE,Q)$  induces a one
parameter subgroup
with $\mu(z,\lambda)=0$ (where $z$ is the corresponding point on
$Z^{ss}$) follows from proposition \ref{giden} and the 
proof of proposition \ref{gmainprop}. 
The fact that the limit point $z_0$
corresponds to $(\CrE_0,Q_0)$ is an easy calculation (see 
\cite[lemma 1.26]{S}). The \tube\ $(\CrE_0,Q_0)$ is semistable by
proposition \ref{conj}.

It is easy to check that $\tilde z$ corresponds to $(\CrE^S,Q^S)$,
and then items (ii) and (iii) follow from the fact that $\tilde
z$ is in $B$.

\end{proof}

\begin{proposition}
\label{conj}
Let $\lambda$ be a 1-PS of $\slv$. 
Let $\slv$ act on $Z$. Asume this action is linearized with respect to
an ample line bundle $\CrH$.
Let $z\in Z^{ss}$. Let $z_0=\lim_{t\to 0}
\lambda (t) \cdot z$. If $\mu(z,\lambda)=0$, then $z_0\in Z^{ss}$.
\end{proposition}

\begin{proof}
This proof was given to us by A. King. We can assume, without loss of 
generality, that the polarization $\CrH$ of $Z$ is very ample, and then $Z$ 
embedds in $\PP (H^0(\CrH)^\vee)$ and $\slv$ acts on $H^0(\CrH)^\vee$.
A point $x\in Z$ is (semi)stable iff its image in $\PP
(H^0(\CrH)^\vee)$ is (semi)stable, and then we can assume
$(Z,\CrH)=(\PP(\CC^n), \SO (1))$, with $\slv$ acting on $\CC^n$.

Let $\pi:\CC^n\setminus \{0\} \to \PP(\CC^n)$ be the projection. Let 
$z\in \PP(\CC^n)$ be a semistable point and $\lambda$ a 1-PS with
$\mu(z,\lambda) =0$. Let $\tilde z \in \CC^n$ be a point in the fibre
$\pi^{-1}(z)$, and let 
$$
\tilde z_0=\lim_{t\to 0} \lambda (t)\cdot \tilde z.
$$ 
This limit exists and it is not the origin because $\mu(z,\lambda)
=0$. We have $z_0:=\lim_{t\to 0} \lambda (t)\cdot  z= \pi(\tilde z_0)$
(by continuity of $\pi$). Assume that the point $z_0$ is unstable.
Then the closure of the orbit of $\tilde z_0$ contains the origin, but this
closure is included in the closure of the orbit of $\tilde z$, and this
doesn't contain the origin because $z$ is semistable. Then $z_0$ is
semistable. Furthermore, $z_0$ cannot be stable because
$\mu(z_0,\lambda) =\mu(z,\lambda)=0$, then $z_0$ is strictly
semistable.

\end{proof}

\section{Properties of \tubes}
\label{properties}

\subsection{Irreducibility of moduli space}
\hfil
\medskip

First we will show that the semistability and stability of \tubes\ are
open conditions.

\begin{proposition}
Let $(\CrE_T,Q_T,\CrN)$ be a flat family of \tubes\ parametrized by
$T$. The subset $T^{s}$ (resp. $T^{ss}$) corresponding to stable 
(resp. semistable) \tubes\ is open.
\end{proposition}

\begin{proof}
Let $m$ and $l$ be large enough so that $\CrV={p^{}_T}_*
(\CrE_T \otimes p^*_X\SO_X(m))$ is locally free, 
$p^*_T \CrV \otimes p^*_X \SO_X(-m) \to \CrE_T$ 
is a surjection and proposition \ref{giden} holds.

Note that the universal family that was constructed on $Z^{ss}$ in the
proof of theorem I can be extended to the set $Z^{good}$ of ``good'' points. 
Arguing as in the proof of theorem I, there is a finite 
open cover $\{U_i\}_
{i\in I}$ of $T$ and morphisms $f_i:U_i \to Z^{good} \subset Z$. These
morphisms depend on the choices made (the choice of local
trivializations of $\CrV$), but the $\slv$ orbit of $f_i(t)$ are
independent of the choices. In particular, the property of $f_i(t)$
belonging to $Z^{ss}$  only depends on $t$. By proposition
\ref{giden}, $f_i(t)$ lies in $Z^{s}$ (resp. $Z^{ss}$) iff 
the \tube\ $(\CrE_t,Q_t)$ is stable (resp. semistable). Then 
$$
T^{s}=\bigcup_{i\in I} f^{-1}_i(Z^{s})
$$
and the openness of $Z^{s}$ in $Z$ proves that $T^{s}$ is open
(the same argument works for $Z^{ss}$).

\end{proof}

\begin{theorem}
Let $X$ be a Riemann surface. Fix $r$, $d$ and $\tau$. Then 
there exists an integer $l_0$ such that if $\deg \CrL > l_0$, then $\overline
\FM_\tau(r,d,\CrL)$ and $\FM_\tau (r,d,\CrL)$ are
irreducible or empty.
\end{theorem}

\begin{proof}
We will construct a flat family of \tubes\ parametrized by an irreducible
scheme $\wt Y$ with the property that every semistable \tube\ of type
$(r,d,\CrL)$ belongs to the family. Then there is a surjective
morphism $\wt Y^{ss} \to \overline \FM_\tau(r,d,\CrL)$, where $\wt
Y^{ss}$ is the open subset representing semistable points, and this 
proves
that $\overline \FM_\tau(r,d,\CrL)$ is irreducible. Repeating this
with the open subset $\wt Y^{s}$ corresponding to stable points, we prove
that $\FM_\tau(r,d,\CrL)$ is also irreducible.

Let $m$ be large enough so that for any semistable \tube\ $(\CrE,Q)$
in $\overline \FM_\tau(r,d,\CrL)$, the
vector bundle $\CrE(m)$ is generated by global sections (corollary 
\ref{gboundcor1}), and such that 
\begin{equation}
\label{gext1const}
2g-2-d-rm<0.
\end{equation}
Note that $m$ only depends on $X$, $r$, $d$ and $\tau$, but not on
$\CrL$. If we choose $r-1$ generic sections of $\CrE(m)$, we have an 
exact sequence
$$
0 \to \SO_X^{\oplus r-1}(-m) \to \CrE \to \CrM(-m) \to 0
$$
where $\CrM$ is a line bundle of degree $d+rm$. 

By standard methods we
can construct a universal family $\CrF_Y$ of extensions of line
bundles of degree $d+(r-1)m$ by $\SO_X^{\oplus r-1}(-m)$. This will be
parametrized by a scheme $Y$ that has a morphism to $\Pic^{d+rm}(X)$,
and the fibre over a line bundle $\CrM$ is naturally isomorphic to
$\Ext^1 (\CrM(-m),\SO_X^{\oplus r-1}(-m))$. Note that to construct this
family we need that the dimension of this $\Ext^1$ group is constant
when we vary $\CrM$, but this is true thanks to (\ref{gext1const}).
Each point $y\in Y$ corresponds to an extension of the form
$$
0 \to \SO_X^{\oplus r-1}(-m) \to \CrF_y \to \CrM(-m) \to 0.
$$
It follows from the argument in the previous paragraph that all vector
bundles in semistable \tubes\ do occur in this family.

Note that, if $(\CrE,Q)$ is a \tube, $Q$ can be thought of as an 
element of $H^0(\Sym^2 \CrF_y^\vee \otimes  \CrL)$.
Now choose $l_0$ large enough so that for any line bundle $\CrL$ of 
degree $\deg(\CrL)>l_0$ the following holds
$$
H^1(\Sym^2 \CrF_y^\vee \otimes  \CrL)=0
$$
for any $y \in Y$. Then $H^0(\Sym^2 \CrF_y^\vee \otimes  \CrL)$
is constant when we vary $y$, and we can construct a (flat) family of
\tubes\ parametrized by $\wt Y=\VV(\Sym^2 \CrF^\vee \otimes p^*_X
\CrL)$, and every semistable \tube\ of type $(r,d,\CrL)$ belongs to
this family.

\end{proof}

\subsection{Orthogonal bundles}
\label{orthogonalbundles}
\hfil
\medskip

An orthogonal bundle is a vector bundle associated to a principal
bundle with (complex) orthogonal structure group. 
Equivalently, it is a \tube\ $(\CrE,Q)$
with  $\CrL=\SO_X$, such that the bilinear form $Q:\Sym ^2
\CrE \to \SO_X$ induces an isomorphism $Q:\CrE \to \CrE^\vee$.
We will call such a \tube\ a smooth \tube.
In this case the \tube\ gives a smooth conic $\CrC_x$ for each point
$x\in X$. Note that the isomorphism $Q:\CrE \to \CrE^\vee$ induces an
isomorphism $\det Q: \det \CrE \to \det \CrE^\vee$, and then
$\deg(\CrE)=0$ (in fact $(\det(\CrE))^{\otimes 2}=\SO_X$).

There is a notion of stability for orthogonal bundles (see \cite{R}):
a bundle $\CrE$ is orthogonal (semi)stable iff for every proper
isotropic subbundle $\CrF$, $\deg(\CrF)(\leq)0$. The notion of stability
that we have defined for \tubes\ depends in principle on a parameter $\tau$, 
but we will show
that in the case of a smooth \tube, the notion of stability
doesn't depend on the particular value of the parameter. In fact we
will prove that
a smooth \tube\ is $\tau$-(semi)stable iff it is (semi)stable as an
orthogonal bundle.

\begin{lemma}
\label{lemmaorthogonal}
Let $(\CrE,Q)$ be a smooth \tube, and let $\CrF$ be a proper vector
subbundle of $\CrE$. Then

(i) There is an exact sequence
$$
0 \to \CrF^\perp \to \CrE \to \CrF^\vee \to 0,
$$ 
and $\deg(\CrF)=\deg(\CrF^\perp)$.

(ii) If $\CrF$ is isotropic ($c_Q(\CrF)\leq 1$), then $\rk(\CrF)=1$.

(iii) If $\rk(\CrF)=1$, then $c_Q(\CrF)\geq 1$
\end{lemma}

\begin{proof}
(i) Follows from the exact sequence
$$
0 \to \CrF^\perp \to \CrE\isom\CrE^\vee \to \CrF^\vee \to 0,
$$ 
and the fact that $\deg(\CrE)=0$.

(ii) Assume that $\rk(\CrF)=2$. Then, in a basis adapted to $\CrF
\subset \CrE$
$$
Q=\left(
\begin{array}{ccc}
0 & 0 & \cdot \\
0 & 0 & \cdot \\
\cdot & \cdot & \cdot
\end{array}
\right)
$$
and then $\det Q = 0$, contradicting the fact that the \tube\ is
smooth.

(iii) If $c_Q(\CrF)=0$, then 
$$
Q=\left(
\begin{array}{ccc}
0 & 0 & 0 \\
0 & \cdot & \cdot \\
0 & \cdot & \cdot
\end{array}
\right)
$$
and then $\det Q = 0$, again contradicting the fact that the \tube\ is
smooth.

\end{proof}

\begin{proposition}
A smooth \tube\ $(\CrE,Q)$ is $\tau$-semistable iff the vector
bundle $\CrE$ is semistable as an
orthogonal bundle. Furthermore, it is $\tau$-stable iff it is stable
as an orthogonal bundle.
\end{proposition}

\begin{proof}
Let $(\CrE,Q)$ be a smooth $\tau$-semistable \tube. Let $\CrF$ be
an isotropic vector subbundle. By lemma \ref{lemmaorthogonal} (ii),
$\rk(\CrF)=1$. We have $\CrF \subset \CrF^\perp$, 
$\rk(\CrF^\perp)=2$ (by
lemma \ref{lemmaorthogonal} (i)), and we check that  $\CrF \subset \CrF^\perp
\subset \CrE$ is a critical filtration. Then $\deg(\CrF) +
\deg(\CrF^\perp) \leq 0$, but
$\deg(\CrF) = \deg (\CrF^\perp)$ (lemma \ref{lemmaorthogonal} (i)),
and then $\deg(\CrF)\leq 0$, which proves that $\CrE$ is semistable as
an orthogonal bundle. Furthermore, if $(\CrE,Q)$ is
$\tau$-stable, then  $\deg(\CrF) + \deg(\CrF^\perp) < 0$, 
$\deg(\CrF)< 0$ and $\CrE$ is stable as an orthogonal bundle.

Conversely, let $\CrE$ be an orthogonal semistable bundle. Let $\CrF$
be any vector subbundle. Following \cite{R} let $\CrN=\CrF \cap
\CrF^\perp$, and let $\CrN'$ be the vector subbundle generated by
$\CrN$. We have an exact sequence
\begin{equation}
\label{orthoshort}
0 \to \CrN' \to \CrF \oplus \CrF^\perp \to \CrM' \to 0
\end{equation}
where $\CrM'$ is the subbundle of $\CrE$ generated by $\CrF +
\CrF^\perp$. We have $\CrM'=(\CrN')^\perp$.

If $\CrN'=0$, then $\CrE=\CrF\oplus \CrF^\perp$, $c_Q(\CrF)=2$, and
$\deg(\CrF) = 0$ (lemma \ref{lemmaorthogonal} (i)). Then 
$$
\frac{\deg(\CrF) - c_Q(\CrF)\tau}{\rk(\CrF)}=\frac{-2\tau}{\rk(\CrF)}<
\frac{-2\tau}{3}=\frac{\deg(\CrE) - 2\tau}{3}.
$$

If $\CrN'\neq 0$, then $\deg(\CrF)=\deg(\CrN')$ (by lemma
\ref{lemmaorthogonal} (i) and the exact sequence (\ref{orthoshort})),
and then $\deg(\CrF)\leq 0$ (because $\CrE$ is orthogonal semistable
and $\CrN'$ is isotropic). If $\rk(\CrF)=2$, then $c_Q(\CrF)=2$ (by
lemma \ref{lemmaorthogonal} (ii)), and if $\rk(\CrF)=1$, then
$c_Q(\CrF)\geq 1$ (by lemma \ref{lemmaorthogonal} (iii)). In any case
$$
\frac {\deg(\CrF) - c_Q(\CrF)\tau}{\rk(\CrF)} \leq \frac{-
c_Q(\CrF)\tau}{\rk(\CrF)} < \frac{-2\tau}{3}=\frac{\deg(\CrE) - 2\tau}{3}
$$
Now let $\CrE_1 \subset \CrE_2 \subset \CrE$ be a critical
filtration. Then $\CrE_1$ is isotropic, $\CrE_2 = \CrE_1 ^\perp$, and
then 
$$
\deg(\CrE_1)+\deg(\CrE_2)=2\deg(\CrE_1) \leq 0,
$$
because $\CrE_1$ is isotropic and $\CrE$ is orthogonal semistable.
This finishes the proof that $(\CrE,Q)$ is $\tau$-semistable. Furthermore, if
$\CrE$ is orthogonal stable, the last inequality is strict, and
we obtain that $(\CrE,Q)$ is $\tau$-stable.

\end{proof}

\section*{Appendix: Hitchin-Kobayashi correspondence for conic 
bundles}
\centerline{\textbf{(By I. Mundet i Riera)}}
\bigskip

In this appendix I use the result in \cite{Mu} to
relate the notion of stability for conic bundles
to the existence of solutions to a certain PDE. This is similar
to the well known relation between stability of vector bundles and
existence of Hermite--Einstein metrics, or between stability of
holomorphic pairs and solutions to the vortex equations. As usual
in the literature, I call such a relation a Hitchin--Kobayashi
correspondence (see \cite{Mu} and the references therein).

Take a non-degenerate conic bundle $Q:\Sym^2\CrE\to\CrL$ on 
a Riemann surface $X$. Let $E$ be the smooth bundle underlying $\CrE$.
We denote $\overline{\partial}_{\CrE}$ the $\overline{\partial}$
operator on $E$ given by $\CrE$.
Fix a metric (In this appendix {\it metric} will always mean
{\it Hermitian metric}). on $\CrL$ and consider the following
equation on a metric $h$ on $E$:
\begin{equation}
\imag\Lambda F_{\overline{\partial}_{\CrE},h}+
\frac{\tau}{2}\frac{Q\otimes Q^{*_h}}{\|Q\|^2_h}=c\Id,
\label{equ0}
\end{equation}
where $F_{\overline{\partial}_{\CrE},h}$ is the curvature of the 
Chern connection of $\overline{\partial}_{\CrE}$ with respect to $h$,
$\Lambda:\Omega^2(X)\to\Omega^0(X)$ is the adjoint of wedging with
the Kaehler form of $X$ and the subscript $h$ in $*$ and $\|\cdot\|$
is to recall that both depend on $h$. Finally, $\tau>0$ and $c$ are
real numbers.

We will take a (rather standard) point of view putting equation 
(\ref{equ0}) inside the setting considered in \cite{Mu}.
Then we will study the existence criterion to solutions of the
equation given in \cite{Mu} applied to this particular case, thus
arriving again at the notion of stability for conic bundles.

Fix a metric $h_0$ in $E$. Let $\SG^c$ be the complex gauge group of
$E$, i.e., the group of its smooth automorphisms covering the identity
on $X$. The group $\SG^c$ acts on the space of
$\overline{\partial}$ operators on $E$ by pullback. So $g\in\SG^c$
sends $\overline{\partial}_{\CrE}$ to $g^*\overline{\partial}_{\CrE}=
g\circ \overline{\partial}_{\CrE}\circ g^{-1}$.
Any metric $h$ on $\CrE$ is the pullback $h=g^*h_0$
by some $g\in\SG^c$. Furthermore, for any metric $h$ and gauge
transformation $g$
$$g\left(F_{(g^{-1})^*\overline{\partial}_{\CrE},h}\right) g^{-1}
=F_{\overline{\partial}_{\CrE},g^*h}
\qquad\mbox{and}\qquad
g\left(\frac{g^{-1}Q\otimes (g^{-1}Q)^{*_h}}{\|g^{-1}Q\|^2_h}
\right) g^{-1}=\frac{Q\otimes Q^{*_{(g^*h)}}}{\|Q\|^2_{g^*h}}.$$
So if $h=g^*h_0$ solves (\ref{equ0}) then, conjugating by $g$, we get 
\begin{equation}
\imag\Lambda F_{(g^{-1})^*\overline{\partial}_{\CrE},h_0}+
\frac{\tau}{2}
\frac{g^{-1}Q\otimes g^{-1}Q^{*_{h_0}}}{\|g^{-1}Q\|^2_{h_0}}=c\Id.
\label{equ1}
\end{equation}
We will now see that equation (\ref{equ1}) on $g\in\SG^c$ is a 
particular case of the equations considered in \cite{Mu}. 

Let $F=\PP(\Sym^2(\CC^3)^*)$.
Take on $F$ the symplectic structure $\tau\omega$, where $\omega$
is the symplectic structure on $F$ obtained from the canonical metric
on $\Sym^2(\CC^3)^*$. 
(For that we view $F$ as the symplectic quotient $F=\mu_0^{-1}(-\imag)/S^1$,
where $\mu_0(z)=-\imag|z|^2$ is the moment map of the action of $S^1$ on
$\Sym^2(\CC^3)^*$.)
Consider on $F$ the action of $U(3,\CC)$ induced by the canonical
action on $\CC^3$. This action is Hamiltonian, and the moment map 
evaluated at $x\in F$ is
\begin{equation}
\mu(x)=-\imag\frac{\tau}{2}\left(\frac{\hat{x}\otimes\hat{x}^*}{\|\hat{x}\|^2}
\right),
\label{mom}
\end{equation}
where $\hat{x}\in\Sym^2(\CC^3)^*$ is any lift of $x$. 

Let $P$ be the $U(3,\CC)$ principal bundle of $h_0$-unitary
frames of $E$, and let $\SF=P\times_{U(3,\CC)}F$. The conic bundle
$Q$ gives a section $\Phi\in\Gamma(\SF)$, and by formula (\ref{mom})
the term in (\ref{equ1}) involving $Q$ is $\imag\mu(\Phi)$. 
Let $\SA$ be the set of
connections on $P$. Let $A=A_{\overline{\partial}_{\CrE},h_0}$ be the
Chern connection. The action of $\SG^c$ on $\SA$ considered in 
\cite{Mu} is as follows: $g\in\SG^c$ sends $A$ to 
$g(A)=A_{g^*\overline{\partial}_{\CrE},h_0}$. 
Finally, since the conic bundle $Q:\Sym^2\CrE\to\CrL$ is
non-degenerate, the pair $(A,\Phi)$ is simple. 
So by the theorem in \cite{Mu} there is a solution $g\in\SG^c$ to equation 
(\ref{equ1}) if and only if $(A,\Phi)$ is $c$-stable.
Furthermore, the metric $g^*h_0$ is unique.

The previous discussion applies also to bundles of quadrics
on projective bundles of arbitrary dimension. In the next section
we will study the $c$-stability condition on any rank and in the next
one we will give a more precise description of $c$-stability 
for conic bundles.

\subsection*{Stability for bundles of quadrics}
\hfil
\medskip

We will suppose from now on
that $\Vol(X)=1$. The stability condition stated in \cite{Mu} 
refers to reductions of the structure group of our bundle to parabolic 
subgroups plus antidominant characters of those parabolic subgroups.
In our case the structure group is $GL(n,\CC)$, so a parabolic reduction
is equivalent to a filtration by subbundles:
$$0\subset E_1\subset\dots\subset E_r=E,$$
where the ranks strictly increase.
The action on $E$ of any antidominant character for this reduction 
is given by a matrix of this form (written using any splitting
$E=E_1\oplus E_2/E_1\oplus\dots\oplus E_r/E_{r-1}$)
\begin{equation}
\chi=\left(
\begin{array}{cccc}
z+m_1+\dots+m_{r-1} & 0 & \dots & 0 \\
0 & z+m_2+\dots+m_{r-1} & \dots & 0 \\
\vdots & \vdots & \ddots & \vdots \\
0 & 0 & \dots & z
\end{array}
\right)-\sum_{k=1}^{r-1} m_k\frac{\rk(E_k)}{\rk(E)}\Id,
\label{caracter}
\end{equation}
where $z$ is any real number and $m_j\leq 0$ are negative real numbers
(strictly speaking this is the action of $\imag$ times an antidominant
character; however, we will ignore this in the sequel. Using the
notation of \cite{Mu} this matrix is $\imag g_{\sigma,\chi}$, where
$\chi$ is an antidominant character).

\textbf{Stability of a quadric.}
To write the stability notion for $(A,\Phi)$ we need to compute
the maximal weight of the action of $\chi$ on the section $\Phi$. 
So fix a point $x\in X$ and write $$0\subset W_1\subset\dots\subset W_r=W$$
the induced filtration in the fibre $W=E_x$ over $x$.
Take a basis $e_1,\dots,e_n$ of $W$ such that for any $1\leq k\leq r$,
$\{e_1,\dots,e_{\rk(W_k)}\}$ is also a basis of $W_k$.
Write $e_1^*,\dots,e_n^*$ the dual basis, so that $Q$ gives
on $W$ the quadratic form
$$Q=Q(x)=\sum_{i\leq j} \alpha_{ij}(e_i^* e_j^*).$$
The action of $\chi$ on $\Sym^2W^*$ diagonalizes in the basis
$\{e_i^* e_j^*\}_{i\leq j}$, and one has
\begin{eqnarray*}
\lefteqn{\chi(e_i^* e_j^*)=}\\
 & & \left(-(2z+2m_I+\dots+2m_{J-1}+m_J+\dots+m_{r-1})
+2\sum_{k=1}^{r-1} m_k\frac{\dim(W_k)}{\dim(W)}\right) (e_i^* e_j^*).
\end{eqnarray*}
Here and in the sequel we follow this convention: the index $I$
(resp. $J$) is the minimum one such that $e_i$ belongs to $W_I$ (resp. $e_j$
belongs to $W_J$). From this one deduces that 
$$\mu(Q(x);\chi)=\max_{\alpha_{ij}\neq 0}
\{-(2z+2m_I+\dots+2m_{J-1}+m_J+\dots+m_{r-1})\}
+2\sum_{k=1}^{r-1} m_k\frac{\dim(W_k)}{\dim(W)}.$$

Define $M_I=-(m_I+\dots+m_{r-1})$. Given two subspaces $W',W''\subset W$, 
we will write $Q(W',W'')=0$ if for any $w'\in W'$ and 
$w''\in W''$, $Q(w',w'')=0$. Otherwise we will write
$Q(W',W'')\neq0$. Then,
\begin{equation}
\mu(Q(x);\chi)=\max_{Q(W_I,W_J)\neq 0}\{M_I+M_J-2z\}
+2\sum_{k=1}^{r-1} m_k\frac{\dim(W_k)}{\dim(W)}\Id.
\end{equation}

\textbf{Stability for the bundle of quadrics.}
The pair $(A,\Phi)$ is by definition
$c$-stable if for any filtration $\sigma$ of $E$ by subbundles
$$0\subset E_1\subset\dots\subset E_r=E$$
and any antidominant character $\chi$ as in (\ref{caracter}) one has
\begin{equation}
\deg(\sigma,\chi)+\tau\int_{x\in X}\mu(Q(x);\chi)-\la\chi,c\Id\ra>0.
\label{estable}
\end{equation}
Here the degree of the pair $(\sigma,\chi)$ is
$$\deg(\sigma,\chi)=z\deg(E)+\sum_{j=1}^{r-1}m_j
\left(\deg(E_j)-\frac{\rk(E_j)}{\rk(E)}\deg(E)\right),$$
and, on the other hand, $\la\chi,c\Id\ra=zc\rk(E)=zcn$. 
The map $Q$ is holomorphic, and the function $\mu(Q(x);\chi)$ is lower
semicontinuous and takes a finite number of values as $x$ moves on $X$. 
Hence, $\mu(Q(x);\chi)$ takes its maximal value in a Zariski open
dense subset of $X$, and so 
$$\int_{x\in X}\mu(Q(x);\chi)=
\Vol(X)\max_{x\in X}\mu(Q(x);\chi)=\max_{x\in X}\mu(Q(x);\chi).$$
For any pair of subbundles $E',E''\subset E$, 
define $Q(E',E'')=\max_{x\in X} Q(E'_x,E''_x).$ Then
\begin{equation}
\max_{x\in X}\mu(Q(x);\chi)=\max_{Q(E_i,E_j)\neq 0}\{M_i+M_j-2z\}+
\sum_{k=1}^{r-1}2m_k\frac{\rk(E_k)}{\rk(E)}.
\label{maxim2}
\end{equation}
Putting everything together (\ref{estable}) becomes
\begin{align}
0 &< z\deg(E)-2\tau z-zc\rk(E)
+\sum_{k=1}^{r-1}m_k
\left(\deg(E_k)-\frac{\rk(E_k)}{\rk(E)}\deg(E)
+2\tau \frac{\rk(E_k)}{\rk(E)}
\right) \notag \\
& +\tau\max_{Q(E_i,E_j)\neq 0}\{M_i+M_j\}.\notag
\end{align}
This must be true for any real number $z$, so the pair $(A,\Phi)$
can only be stable if 
$$c=\frac{\deg(E)-2\tau }{\rk(E)}.$$
Define now 
$$d_k=\deg(E_k)-\frac{\rk(E_k)}{\rk(E)}\deg(E)
+2\tau \frac{\rk(E_k)}{\rk(E)}.$$
Then the stability condition reduces to
\begin{equation}
\sum_{k=1}^{r-1} m_kd_k+\tau\max_{Q(E_i,E_j)\neq 0}\{M_i+M_j\}>0.
\label{defi}
\end{equation}
And this must hold for any choice of (not all zero) negative numbers
$m_1,\dots,m_{r-1}$.

\section{The case $\rk(E)=3$}
In the sequel we will use the following notation. If $E'$ is a
vector bundle and $\alpha$ is any real number,
$$\mu_{\alpha}(E')=\frac{\deg(E')-\alpha}{\rk(E')}.$$

In this section we assume that $\rk(E)=3$. Hence, $Q$ describes a 
bundle of conics in a bundle of projective planes $\PP(E)$ on $X$. 
Recall that we assume that $Q$ is (generically) non-degenerate.
We have seen above that the pair $(A,\Phi)$ cannot be $c$-stable unless
$$c=\mu_{2\tau}(E).$$
Suppose this holds. Now, according to formula (\ref{defi}), 
$(A,\Phi)$ is stable if and only if for any filtration 
$0\subset E_1\subset E_2\subset E$ and
for any pair of (not all zero) real numbers $m_1,m_2\leq 0$,
\begin{equation}
m_1d_1+m_2d_2+\tau\max_{Q(E_i,E_j)\neq 0}\{M_i+M_j\}>0,
\end{equation}
where, as before, 
$d_k=\deg(E_k)-\frac{\rk(E_k)}{\rk(E)}\deg(E)+2\tau \frac{\rk(E_k)}{\rk(E)}$.
There are three cases to consider:
\begin{itemize}

\item $Q(E_1,E_1)=Q(E_1,E_2)=0$, $Q(E_2,E_2)\neq 0$.
Geometrically, $E_1$ gives fibrewise a point on the conic
and $E_2$ a tangent line to the conic at the point given by $E_1$.
In this case, 
$$\max_{Q(E_i,E_j)\neq 0}\{M_i+M_j\}=\max\{-2m_2,-m_1,-m_2\}.$$
Hence,
\begin{align*}
0 &> d_1-\tau=\deg(E_1)-\frac{\rk(E_1)}{\rk(E)}\deg(E)+
2\tau\frac{\rk(E_1)}{\rk(E)}-\tau, \\
0 &> d_2-2\tau=\deg(E_2)-\frac{\rk(E_2)}{\rk(E)}\deg(E)+
2\tau\frac{\rk(E_2)}{\rk(E)}-2\tau, \\
0 &> d_1+d_2-2\tau.
\end{align*}
Simplifying, we obtain the following conditions:
$$\mu_{\tau}(E_1)<\mu_{2\tau}(E)
\mbox{, }\qquad\mu_{2\tau}(E_2)<\mu_{2\tau}(E)
\mbox{, }\qquad\deg(E_1)+\deg(E_2)<\deg(E).$$
\item 
$Q(E_1,E_1)=0$, $Q(E_1,E_2)\neq 0$ ($\Rightarrow\ Q(E_2,E_2)\neq 0$).
Geometrically, $E_1$ is a point on the conic and $E_2$ a line
passing through $E_1$ but generically not tangent to the conic.
In this case, 
$$\max_{Q(E_i,E_j)\neq 0}\{M_i+M_j\}=-m_1-2m_2.$$    
Hence,
\begin{align}
0 &> d_1-\tau=\deg(E_1)-\frac{\rk(E_1)}{\rk(E)}\deg(E)+
2\tau\frac{\rk(E_1)}{\rk(E)}-\tau, \notag \\
0 &> d_2-2\tau=\deg(E_2)-\frac{\rk(E_2)}{\rk(E)}\deg(E)+
2\tau\frac{\rk(E_2)}{\rk(E)}-2\tau.\notag
\end{align}
Simplifying, we obtain the following two conditions:
$$\mu_{\tau}(E_1)<\mu_{2\tau}(E)
\mbox{, }\qquad\mu_{2\tau}(E_2)<\mu_{2\tau}(E).$$
\item 
$Q(E_1,E_1)\neq 0$ ($\Rightarrow\ Q(E_1,E_2)\neq 0
\mbox{ and }Q(E_2,E_2)\neq 0$).
Geometrically, $E_1$ gives a point generically not on the conic
and $E_2$ any line through $E_1$.
In this case, 
$$\max_{Q(E_i,E_j)\neq 0}\{M_i+M_j\}=-2m_1-2m_2.$$
Hence,
\begin{align}
0 &> d_1-2\tau=\deg(E_1)-\frac{\rk(E_1)}{\rk(E)}\deg(E)+
2\tau\frac{\rk(E_1)}{\rk(E)}-2\tau, \notag \\
0 &> d_2-2\tau=\deg(E_2)-\frac{\rk(E_2)}{\rk(E)}\deg(E)+
2\tau\frac{\rk(E_2)}{\rk(E)}-2\tau. \notag
\end{align}
Simplifying, we obtain the following two conditions:
$$\mu_{2\tau}(E_1)<\mu_{2\tau}(E)
\mbox{, }\qquad\mu_{2\tau}(E_2)<\mu_{2\tau}(E).$$

\end{itemize}

\vskip 1cm

In conclusion, and as claimed at the beginning, the condition of 
$c$-stability obtained from studying equation (\ref{equ0}) coincides 
with that of stability obtained from the GIT construction of the
moduli space of conic bundles.

\bigskip
\textbf{Acknowledgements.}
We would like to thank R. Hern\'andez, A. King, S. Ramanan and 
C.S. Seshadri for discussions and comments.

\bigskip

T. G\'omez

\texttt{tomas@math.tifr.res.in}

School of Mathematics

Tata Institute of Fundamental Research

Homi Bhabha Road, 400 005 Mumbai (India)

\bigskip

I. Sols

\texttt{sols@eucmax.sim.ucm.es}

Departamento de Algebra

Facultad de Ciencias Matem\'aticas

Universidad Complutense de Madrid, 28040 Madrid (Spain)

\bigskip

I. Mundet i Riera

\texttt{ignasi.mundet@uam.es}

Dep. Matem\'aticas, Facultad de Ciencias

Universidad Aut\'onoma de Madrid, Madrid (Spain)

\end{document}